[1994/12/01]
\documentclass[11pt]{amsart}
\chardef\bslash=`\\ 

\usepackage{enumerate}

\usepackage{vaucanson-g}
\usepackage{pstricks}
\usepackage{verbatim}

\newcommand{\NN}{\mathbb{N}}
\newcommand{\ZZ}{\mathbb{Z}}

\newcommand{\AAA}{\mathcal{A}}
\newcommand{\BBB}{\mathcal{B}}
\newcommand{\CCC}{\mathcal{C}}
\newcommand{\LLL}{\mathcal{L}}
\newcommand{\DDD}{\mathcal{D}}

\newcommand{\GGG}{\mathcal{G}}
\newcommand{\htop}{h_\mathrm{top}}
\newcommand{\ulim}{\varlimsup}

\newcommand{\eps}{\varepsilon}

\newcommand{\mc}{\mathcal}
\newcommand{\symdiff}{\bigtriangleup}
\newtheorem{definition}{Definition}[section]
\newtheorem{theorem}{Theorem}[section]
\newtheorem{lemma}[theorem]{Lemma}
\newtheorem{proposition}[theorem]{Proposition}
\newtheorem{corollary}[theorem]{Corollary}
\newtheorem{thma}{Theorem}
\newtheorem*{thma*}{Theorem}

\theoremstyle{remark}
\newtheorem*{remark}{Remark}
\newtheorem*{problem}{Problem}

\newcommand{\V}{(W)}

\newcommand{\F}{(S)}
\newcommand{\Fp}{(Per)}

\DeclareMathOperator{\Per}{Per}

\numberwithin{equation}{section}

\newcommand{\eval}[2][\right]{\relax
  \ifx#1\right\relax \left.\fi#2#1\rvert}

\begin{document}
\title[Intrinsic ergodicity beyond specification]{Intrinsic ergodicity beyond specification:  $\beta$-shifts, $S$-gap shifts, and their factors}
\author[V. Climenhaga]{Vaughn Climenhaga}
\address{Department of Mathematics\\ University of Maryland\\College Park, MD 20742}
\email{climenhaga@math.umd.edu}
\author[D.J. Thompson]{Daniel J. Thompson}
\address{Department of Mathematics\\ Pennsylvania State University\\ University Park, PA 16802}
\email{thompson@math.psu.edu}

\date{\today}
\begin{abstract}
We give  sufficient conditions for a shift space $(\Sigma,\sigma)$ to be intrinsically ergodic, along with sufficient conditions for every subshift factor of $\Sigma$ to be intrinsically ergodic.  As an application, we show that every subshift factor of a $\beta$-shift is intrinsically ergodic, which answers an open question included in Mike Boyle's article ``Open problems in symbolic dynamics''.  We obtain the same result for $S$-gap shifts, and describe an application of our conditions to more general coded systems.  One novelty of our approach is the introduction of a new version of the specification property that is well adapted to the study of symbolic spaces with a non-uniform structure.
\end{abstract}
\maketitle

\section{Introduction}
We study uniqueness of a measure of maximal entropy (or \emph{mme} for short) in the context of symbolic dynamics.  Dynamical systems with a unique mme are called \emph{intrinsically ergodic}. Determining which dynamical systems are intrinsically ergodic is a central problem at the interface of ergodic theory and topological dynamics \cite{Ho2, Ho3, Gu72, We, Bu2, BG, BuFi, Bo7}, and is a prototypical result for the thermodynamic formalism, a powerful tool for studying the statistical properties of a system.

Parry~\cite{wP64} and Weiss~\cite{We, bW73} established intrinsic ergodicity for topologically transitive shifts of finite type (SFTs), and all their subshift factors (sofic shifts).  Given a property defining a class of shifts, it is natural to ask whether this property implies intrinsic ergodicity, and whether it is preserved by passing to factors.  This is known to be the case for the specification property~\cite{Bo7}, but beyond specification, less is known.

Our motivating examples are the family of $\beta$-shifts and the family of $S$-gap shifts, which do not have specification (see \S \ref{beta}).  In particular, we answer the following open problem posed by Klaus Thomsen, which is Problem 28.2 of Mike Boyle's article ``Open problems in symbolic dynamics''~\cite{Boy}:
\begin{problem}
Must a subshift factor of a $\beta$-shift be intrinsically ergodic?
\end{problem}
\begin{thma}\label{thm:beta-factor}
Every subshift factor of a $\beta$-shift or $S$-gap shift is intrinsically ergodic. Moreover, the unique measure of maximal entropy can be characterised as the (well defined) weak* limit as $n \rightarrow \infty$ of $\delta$-measures evenly distributed across periodic points of period at most $n$.
\end{thma}

Theorem~\ref{thm:beta-factor} is proved via our more general main result (Theorem~\ref{thm:main}), which provides natural sufficient conditions for an abstract shift space to be intrinsically ergodic.  These conditions may be viewed as a weakening of the specification property, and are well behaved under the operation of taking factors.  

Our conditions take a particularly simple form for the class of coded systems, which includes $\beta$-shifts and $S$-gap shifts.  A shift space $X\subset \Sigma_p = \{1,\dots,p\}^\ZZ$ is \emph{coded} if there exists a countable collection of finite words, called generators, such that $X$ is the closure of the set of sequences obtained by freely concatenating the generators.  Given a set of generators for a coded system, let $c_n$ be the number of words 
of length $n$ that appear either at the beginning or the end
of some generator. 
\begin{thma}\label{thm:coded}
Let $(X,\sigma)$ be a coded shift and let $c_n$ be as above. 
\begin{enumerate}
\item If $\ulim_{n\to\infty} \frac 1n \log c_n < \htop(X,\sigma)$, then $(X,\sigma)$ is intrinsically ergodic.
\item If $\lim_{n\to\infty} \frac 1n \log c_n = 0$, then every subshift factor of $(X,\sigma)$ is intrinsically ergodic.
\end{enumerate}
Moreover, under these conditions, the unique measure of maximal entropy is the weak* limit of $\delta$-measures evenly distributed on periodic points of period at most $n$.
\end{thma}

A crucial ingredient in our approach, and a result of interest in its own right, is that under our conditions, the unique measure of maximal entropy satisfies a certain Gibbs property. We also give a sufficient condition for the unique measure of maximal entropy to be obtained as the weak limit of periodic orbit measures. 

We note that by expansivity, a subshift over a finite alphabet always has at least one measure of maximal entropy, so the main content of the theorem is uniqueness.  There are many examples of non-intrinsically ergodic subshifts in the literature \cite{De, kP86}. One can even construct minimal subshifts with arbitrarily many measures of maximal entropy \cite[Example 27.2]{De}.



Beyond the shifts with specification, various classes of shift spaces have been shown to be intrinsically ergodic, but none of these classes are closed under passing to factors.  For example, a class of shifts called \emph{almost sofic} was introduced by Petersen in~\cite{kP86}.  Many of these shifts are intrinsically ergodic, but not all, and Petersen gives an explicit example of an almost sofic shift that is intrinsically ergodic but which nevertheless has a subshift factor with more than one measure of maximal entropy.

Thus the class of shifts that are both almost sofic and intrinsically ergodic is not closed under factors, and Petersen observed that ``a useful class of almost sofic, intrinsically ergodic systems which contains the sofic systems and is closed under the usual dynamical operations such as passing to factors has not yet been identified.''  One merit of our approach is that we are able to describe a class of intrinsically ergodic systems which meets Petersen's criterion of being closed under passing to factors.  

We compare our approach with existing techniques for proving intrinsic ergodicity beyond specification, focusing on the $\beta$-shifts, although we emphasise that our techniques and results work in a more general setting.  Each value of $\beta>1$ determines a subshift $\Sigma_\beta$. The $\beta$-shifts are a very natural and explicit family of shift spaces which, for generic values of $\beta$, are not contained in the usual classes of shift spaces where standard techniques apply  (eg.\ SFTs, sofic shifts, shifts with specification).  Interest in the $\beta$-shift arises from its role as the coding space for the $\beta$-transformation, which has a deep connection with number theory (see \S \ref{beta}).

Intrinsic ergodicity for arbitrary $\beta$-shifts was established independently by Hofbauer~\cite{Ho2} and by Walters~\cite{Wa2}.  Hofbauer's approach relies on modeling the $\beta$-shift by a countable state topological Markov chain with strong recurrence properties, and using a version of the Perron--Frobenius theorem to establish intrinsic ergodicity.  The approach taken by Walters, on the other hand, applies transfer operator methods (related to the Perron--Frobenius approach above) directly to the $\beta$-shift.

Both of these approaches have been extensively generalised, and have proved very successful in a multitude of situations. However, it seems problematic to adapt these methods to the operation of taking factors.  We develop another approach, which does not use Perron--Frobenius theory, and which has more in common with Bowen's elegant proof that expansive maps with specification are intrinsically ergodic~\cite{Bo7}.

Given a shift space $X$, we write $\LLL$ for the language of $X$---that is, the collection of all finite words that appear in sequences $x\in X$.  In this context, specification is the ability to use a connecting word of a fixed length to glue together any two words $v, w$ from the language of the shift space---that is, the existence of $t \in \NN$ such that given any $v,w\in \LLL$, there is a word $x\in \LLL$ with length $|x| = t$ for which the concatenation $vxw$ is once again in $\LLL$.  For shifts with specification, Bowen's proof proceeds by using combinatorial arguments to establish a Gibbs property for a certain measure of maximal entropy, and then using this Gibbs property to prove uniqueness.

For more general shifts, topological transitivity guarantees the existence of some  $x\in \LLL$ so that $vxw \in \LLL$, but $x$ may be arbitrarily long; this is the case for generic $\beta$-shifts.  This necessitates a new approach to the estimates in Bowen's proof, which no longer hold in their original form.  We overcome this difficulty by considering a collection of ``good'' words $\GGG \subset \LLL$ on which specification holds.  We use the structure of $\Sigma_\beta$ to characterise the obstructions that prevent a word from being good: the words that do not belong to $\GGG$ are precisely those that end in a word taken from a certain smaller collection of words $\CCC^s$.  We are able to describe  the  collection  $\CCC^s$ very explicitly, and the growth rate of the number of words of length $n$ in $\CCC^s$ is subexponential.  This allows us to prove that a uniformly positive proportion of all the words of length $n$ in the language are `good', which in turn allows us to establish a weakened Gibbs property and prove uniqueness.

For our more general results, including $S$-gap shifts and general coded shifts, we pursue a similar strategy.  Given a collection of words $\GGG$ that satisfies the specification property, 
we characterise the words that do not belong to $\GGG$ explicitly by a collection of their possible endings $\CCC^s$ (suffixes) and a collection of their possible beginnings $\CCC^p$ (prefixes). We require that every word can be extended to a word in $\GGG$ in a suitably uniform manner. If the number of words of length $n$ arising from the collections $\CCC^s$ and $\CCC^p$ grows more slowly than the topological entropy, then the shift is intrinsically ergodic.


We mention that our work is in a similar spirit to Buzzi's work on shifts of quasi-finite type \cite{Bu3}, of which the $\beta$-shifts are a prime example. The key property he assumes is that the number of words that are ``constraints'' grows slower than the entropy.  This property is philosopically similar to our requirement that a certain collection of words grows slower than the entropy.  Nevertheless, our conditions do not seem to imply the quasi-finite type property, and definitely do not follow from it.  Buzzi was able to show (among many other things) that q.f.t.\ shifts have finitely many ergodic measures of maximal entropy, and gave counter-examples to uniqueness.

Asking for the specification property to hold only for words taken from a suitable proper subset of $\LLL$ is a key innovation in our approach, and provides the necessary flexibility to deal with shift spaces whose behaviour is a long way from being Markov. There has been a resurgence in interest in specification properties recently, due to important contributions by Pfister and Sullivan (almost specification \cite{PfS, Ya, Tho7}), and Varandas (non-uniform specification \cite{Va}). Our work is very much in the spirit of these developments, giving another direction in which to weaken the specification property in order to apply to a wider range of examples.

In \S \ref{results}, we collect our definitions and state our main result, a condition for intrinsic ergodicity, together with results on how this condition behaves under factors.  In \S \ref{beta}, we discuss in detail the application to $\beta$-shifts, $S$-gap shifts, and their factors, showing that Theorem~\ref{thm:beta-factor} follows from our main results.  In \S \ref{examples}, we discuss coded systems and derive Theorem~\ref{thm:coded}.  In \S \S \ref{proof}-\ref{factors}, we prove our main results on intrinsic ergodicity, and their behaviour under factors.

\section{Definitions and statement of result}\label{results}

A topological dynamical system is a compact metric space $X$ together with a continuous map $f\colon X\to X$. Let $\mathcal M_f(X)$ denote the space of  $f$-invariant probability measures on $X$.  We write $\htop(X, f)$ for the topological entropy of the dynamical system, and $h_\mu(f)$ for the measure-theoretic entropy of $\mu \in \mathcal M_f(X)$. The variational principle \cite[Theorem 8.6]{Wa} states that
\[
\htop(X,f) = \sup \{ h_\mu(f) \mid \mu \in \mathcal M_f(X)\}.
\]
An invariant probability measure that attains this supremum is called a \emph{measure of maximal entropy} (or \emph{mme} for short).  If such a measure exists and is unique, the system is called \emph{intrinsically ergodic}.

\subsection{Languages for shifts}

We begin by recalling the relationship between shift spaces and languages.  We refer the reader to ~\cite{BH,LM} for further background and proofs.

Fix an integer $p\geq 2$ and let $\{1,\dots,p\}^{< \NN}$ be the collection of all finite words in the symbols $1, \dots, p$.  Juxtaposition denotes concatenation---that is, given two words $v=v_1\cdots v_m$ and $w=w_1\cdots w_n$, we write $vw = v_1 \cdots v_m w_1 \cdots w_n$.

A \emph{one-sided language} $\LLL \subset \{1, \dots, p\}^{< \NN}$ is a collection of words such that
\begin{enumerate}
\item if $w\in \LLL$ and $v$ is a subword of $w$, then $v\in \LLL$;
\item if $w\in \LLL$, then there exists $a\in \{1,\dots,p\}$ such that $wa\in \LLL$.
\end{enumerate}
Given a one-sided language $\LLL$, let $X=X_\LLL \subset \Sigma_p^+$ be the collection of all sequences $x_1 x_2 \ldots$ such that
\begin{equation}\label{eqn:realword}
x_i x_{i+1} \dots x_{j-1} x_j \in \LLL
\end{equation}
for every $1\leq i\leq j < \infty$.  Then $X$ is a closed $\sigma$-invariant set, where $\sigma$ is the usual shift operator defined on $\{1,\dots,p\}^{\NN}$; this is the one-sided shift associated with the language $\LLL$.

The construction also runs in the converse direction:  given a one-sided shift $X\subset \Sigma_p^+$, the set of all words that appear in sequences $x\in X$ is a one-sided language.  This gives a one-to-one correspondence between one-sided shift spaces and one-sided languages.

If $\LLL$ is a one-sided language that satisfies the additional condition
\begin{enumerate}
\setcounter{enumi}{2}
\item if $w\in \LLL$, then there exists $a\in \{1,\dots,p\}$ such that $aw\in \LLL$,\label{two-side}
\end{enumerate}
then we say that $\LLL$ is a \emph{two-sided language}.  Let $\hat X = \hat X_\LLL \subset \Sigma_p$ be the collection of all doubly infinite sequences $\ldots x_{-1} x_0 x_1 \ldots$ such that~\eqref{eqn:realword} holds for every $-\infty < i \leq j < \infty$.  Then $\hat X_\LLL$ is a closed $\sigma$-invariant set, where $\sigma$ is the usual shift operator defined on $\{1,\dots,p\}^{\ZZ}$; this is the two-sided shift associated with the language $\LLL$, and is the natural extension of $X_\LLL$.  As with one-sided shifts, the correspondence runs both ways. The following well known proposition (whose proof we give in Section~\ref{sec:measures} for completeness) shows that intrinsic ergodicity for $X_\LLL$ is equivalent to intrinsic ergodicity for $\hat X_\LLL$.

\begin{proposition} \label{prop:measures}
The invariant measures of $X_\LLL$ and $\hat X_{\LLL}$ can be identified by a natural entropy preserving bijection.
\end{proposition}
Given a one-sided language $\LLL$, let $\hat\LLL$ be the union of all subsets of $\LLL$ that satisfy (\ref{two-side}).  Then $X_{\hat\LLL} = \bigcap_{n\geq 0} \sigma^n(X_\LLL)$, and so $X_{\hat\LLL}$ and $X_\LLL$ have the same space of invariant measures.  Thus it suffices to consider two-sided languages.

Let $|w|$ denote the length of a word $w$, and denote by $\LLL_n$ the collection of all words of length $n$ in $\LLL$.  For one-sided shifts, there is a one-to-one correspondence between words $w \in \LLL$ and central cylinders
\[
[w] := \{x\in X_\LLL \mid x_i = w_i \text{ for all } 1\leq i\leq |w| \}.
\]
We have a similar correspondence between words and cylinders for two-sided shifts, with the caveat that we must keep track of where the cylinder begins: given $w\in \LLL$ and $k\in \ZZ$, we define the cylinder
\[
{}_k[w] := \{x \in \hat X_\LLL \mid x_{k+i-1} = w_i \text{ for all } 1 \leq i\leq |w| \}.
\]
We define the central cylinder for $w \in \LLL_n$ to be $_k [w]$ where $k = -\lfloor n/2\rfloor$.

Given a collection of words $\DDD \subset \LLL$ and $n\geq 1$, let $\DDD_n = \DDD \cap \LLL_n$ be the set of words of length $n$ in $\DDD$.  We denote the growth rate of the number of words of length $n$ in $\DDD$ by
\begin{equation}\label{eqn:growth}
h(\DDD) = \ulim_{n\to\infty} \frac 1n \log \# \DDD_n.
\end{equation}
The correspondence between words and cylinders implies that
\[
h(\LLL) = \htop(X_\LLL, \sigma).
\]

Given collections of words $\AAA,\BBB \subset \LLL$, we will write
\[
\AAA \BBB = \{ vw \in \LLL \mid v\in \AAA, w\in \BBB \}.
\]
Note that only words in $\LLL$ are included in the concatenation $\AAA\BBB$.  It may be the case that $\AAA$ and $\BBB$ are both non-empty, but a word in $\BBB$ cannot follow a word in $\AAA$, and so $\AAA \BBB = \emptyset$.

We will occasionally write $0^n$ to denote the word that contains the symbol $0$ repeated $n$ times.  Usually, however, superscripts will denote indices---that is, we write $w^1, w^2, \dots$ for a collection of words, so as to reserve the notation $w_i$ for the $i$th entry of the word $w$.

\subsection{Specification properties}

There are a variety of properties in the literature that go by the name ``specification''.  They all have to do with the ability to approximate arbitrary orbit segments by a single trajectory. In the standard definition of specification due to Bowen, the time spent transitioning between orbit segments has a fixed length, independent of the length of the orbit segments. For symbolic spaces, this corresponds to being able to freely concatenate words using connecting words of fixed length. There are a number of variations on this definition, both in the classical and recent literature \cite{Bo7, De, PfS, Va}.

We formulate specification properties that apply only to a subset of the space. Our definition applies to naturally defined subsets of many examples, such as $\beta$-shifts, that do not have specification.

\begin{definition}
Let $\LLL$ be a language (one- or two-sided) and consider a subset $\GGG \subset \LLL$.  Fix $t\in \NN$; either of the following conditions defines a \emph{specification property on $\GGG$} with gap size $t$.
\begin{description}
\item[(S)] For all $m\in \NN$ and $w^1,\dots,w^m\in \GGG$, there exist $v^1,\dots,v^{m-1}\in \LLL$ such that $x := w^1 v^1 w^2 v^2 \cdots v^{m-1} w^m \in \LLL$ and $|v^i| = t$ for all $i$.
\item[(Per)] Condition \F\ holds, and in addition, the cylinder $[x]$ contains a periodic point of period exactly $|x| + t$.
\end{description}
\end{definition}

When $\GGG = \LLL$, we recover the usual specification property of Bowen.

\begin{remark}
We stress that, in our definition,  we only ask that $x \in \LLL$. We do not require that $x \in \GGG$. 
\end{remark}

\begin{remark}
A natural variant on this definition is to allow the connecting words $v^i$ in \F\ to satisfy $|v^i| \leq t$ for all $i$, rather than requiring the equality $|v^i|=t$. We refer to this as \V-specification.  Theorem~\ref{thm:main} and Proposition~\ref{prop:factorcat} below still hold if \F-specification is replaced with \V-specification (one only needs to be a little more careful in the proof of Lemma \ref{GibbsG}), while Proposition~\ref{prop:posent} and Theorem~\ref{thm:CPE} require the stronger condition \F.  In general, \Fp\ is stronger than \F, which is in turn stronger than \V; however, all three are equivalent in the case $t=0$. 
\end{remark}




\subsection{Statement of Results}

We consider languages $\LLL$ admitting a decomposition $\LLL = \CCC^p \GGG \CCC^s$---that is, there are collections of words $\CCC^p, \GGG, \CCC^s \subset \LLL$ such that every word in $\LLL$ can be written as a concatenation of a word from $\CCC^p$ (a prefix), a word from $\GGG$ (a ``good'' core), and a word from $\CCC^s$ (a suffix), in that order (there may be more than one way to do this).

Given such a decomposition, we define collections of words $\GGG(M)$ for each $M\in \NN$ by
\[
\GGG(M) = \{ uvw \mid u \in \CCC^p, v\in \GGG, w\in \CCC^s, |u|\leq M, |w|\leq M \}.
\]
Observe that $\bigcup_{M\geq 1} \GGG(M) = \LLL$.

Let $\Per(n) = \{x\in X \mid \sigma^k(x) = x \text{ for some } 1\leq k\leq n\}$ be the collection of periodic points of period at most $n$, and write
\begin{equation}\label{eqn:permeas}
\mu_n = \frac{1}{\# \Per(n)} \sum_{x \in \Per(n)} \delta_x
\end{equation}
for the probability measures evenly distributed across the points in $\Per(n)$.

The following results apply equally to one-sided and two-sided shift spaces.

\begin{thma}\label{thm:main}
Let $(X,\sigma)$ be a shift space whose language $\LLL$ admits a decomposition $\LLL = \CCC^p \GGG \CCC^s$, and suppose that the following conditions are satisfied:
\begin{enumerate}[(I)]
\item $\GGG$ has \F-specification. \label{glue}
\item $h(\CCC^p \cup \CCC^s) < h(\LLL) = \htop(X,\sigma)$.\label{slow}
\item For every $M \in \NN$, there exists $\tau$ such that given $v \in \GGG(M)$, there exist words $u, w$ with $|u| \leq \tau, |w| \leq \tau$ for which $uvw \in \GGG$. \label{extend}
\end{enumerate}
Then $(X,\sigma)$ is intrinsically ergodic.  If $\GGG$ has \Fp-specification, then the sequence of probability measures~\eqref{eqn:permeas} converges to the unique measure of maximal entropy.
\end{thma}

\begin{remark}
Condition \eqref{extend} says that every word can be extended to a word in $\GGG$, and that the length required to do so is controlled by the prefix length and suffix length. Our proof will show that conditions \eqref{glue} and \eqref{slow} imply $h(\GGG) = h(\LLL)$.
\end{remark}

\begin{remark}
The assumption that $\GGG$ has \Fp-specification guarantees that the periodic orbit measures~\eqref{eqn:permeas} are well defined.  It is possible for a collection of words to satisfy \F-specification without $X$ containing any periodic orbits, in which case the measures in~\eqref{eqn:permeas} are not well defined. 
\end{remark}

The following proposition tells us how the decomposition $\LLL = \CCC^p \GGG \CCC^s$ and its properties behave under factors.

\begin{proposition}\label{prop:factorcat}
Let $(\tilde{X},\tilde{\sigma})$ be a shift factor of $(X,\sigma)$, and denote the corresponding languages by $\tilde\LLL$ and $\LLL$.  If $\LLL$ admits a decomposition $\LLL = \CCC^p \GGG \CCC^s$, then $\tilde\LLL$ admits a decomposition $\tilde\LLL = \tilde\CCC^p \tilde\GGG \tilde\CCC^s$ such that
\begin{enumerate}
\item If $\GGG$ has \F-specification, then $\tilde \GGG$ has \F-specification;
\item
If $\GGG$ has \Fp-specification, then $\tilde \GGG$ has \Fp-specification;
\item $h(\tilde\CCC^p \cup \tilde\CCC^s) \leq h(\CCC^p \cup \CCC^s)$.
\item If $\LLL$ satisfies \eqref{extend}, then so does $\tilde \LLL$.
\end{enumerate}
\end{proposition}

Combining Theorem~\ref{thm:main} with Proposition~\ref{prop:factorcat} gives the following result.

\begin{corollary} \label{cor:factordawg}
Let $(X,\sigma)$ be a shift space whose language admits a decomposition $\CCC^p \GGG \CCC^s$ satisfying (\ref{glue}), (\ref{slow}) and (\ref{extend}), and let $(\tilde X,\sigma)$ be a subshift factor of $(X,\sigma)$  such that $\htop (\tilde X,\sigma) > h(\CCC^p \cup \CCC^s)$. Then $\tilde X$ is intrinsically ergodic.  If $\GGG$ has \Fp-specification, then the sequence of probability measures~\eqref{eqn:permeas} converges to the unique measure of maximal entropy for $(\tilde X,\sigma)$.
\end{corollary}

There are a number of important examples for which the collection of prefixes and suffixes grows subexponentially.  Furthermore, we have the following dichotomy for systems satisfying Conditions~\eqref{glue} and~\eqref{extend}.

\begin{proposition}\label{prop:posent}
Let $(X,\sigma)$ be a shift space whose language has a decomposition $\LLL = \CCC^p \GGG \CCC^s$ satisfying \eqref{glue} and \eqref{extend}.  Then either $X$ has positive entropy or $X$ comprises a single periodic orbit.
\end{proposition}

Since Conditions \eqref{glue} and \eqref{extend} are preserved by factors (Proposition~\ref{prop:factorcat}), we obtain the following result, which gives a broad class of intrinsically ergodic systems that is closed under taking factors.

\begin{thma}\label{thm:CPE}
Let $(X,\sigma)$ be a shift space whose language admits a decomposition $\CCC^p \GGG \CCC^s$ satisfying conditions \eqref{glue}, \eqref{extend} and
\begin{enumerate}[(I$^\prime$)]
\setcounter{enumi}{1}
\item $h(\CCC^p \cup \CCC^s) =0$.\label{zero}
\end{enumerate}
Then every subshift factor of $(X,\sigma)$ is intrinsically ergodic.  If $\GGG$ has \Fp-specification, then the sequence of probability measures~\eqref{eqn:permeas} converges to the unique measure of maximal entropy for $(\tilde X,\sigma)$.
\end{thma}

In Section~\ref{beta}, we will show that every $\beta$-shift and $S$-gap shift has a language with a decomposition satisfying the conditions of Theorem~\ref{thm:CPE}, which proves Theorem~\ref{thm:beta-factor}.

\section{Application to $\beta$-shifts, $S$-gap shifts, and their factors} \label{beta}
\subsection{$\beta$-shifts}

We recall some facts about $\beta$-shifts; further information can be found in \cite{Pa, Joh, PfS, Maia, Thom05, Tho7}, among others.  For any $\beta$-shift, we describe a decomposition of the language which satisfies the hypotheses of Theorem \ref{thm:CPE}, whence the first part of Theorem~\ref{thm:beta-factor} follows.

Fix a real number $\beta >1$, and let $b = \lceil \beta \rceil$ be the smallest integer greater than or equal to $\beta$.  The $\beta$-shift $\Sigma_\beta \subset \{0, \ldots, b\}^\NN$ is the natural symbolic space associated to the $\beta$-transformation $f_\beta: [0,1) \mapsto [0,1)$ given by
\[
f_\beta(x) = \beta x \pmod 1.
\]
There is a uniquely determined sequence $w(\beta) = (w_j(\beta))_{j=1}^\infty$ that is the lexicographic supremum over all solutions to the equation
\[
\sum_{j=1}^{\infty} w_j (\beta) \beta^{-j} = 1.
\]
The $\beta$-shift  can be characterised by
\[
x \in \Sigma_\beta \iff \sigma^k (x) \preceq w(\beta) \text{ for all } k \geq 1,
\]
where $\preceq$ denotes the lexicographic ordering on $\{0, \ldots, b\}^\NN$.  In particular, for every $k$, $\sigma^k (w(\beta)) \preceq w(\beta)$.  

A $\beta$-shift is sofic if and only if $\beta$ is eventually periodic, and has specification if and only if $w(\beta)$ does not contain arbitrarily long strings of zeroes; the set of $\beta$ with this property has Lebesgue measure zero.  Generically, then, a $\beta$-shift is not sofic and does not possess the specification property \cite{BM, Schm}.

\begin{figure}[htbp]
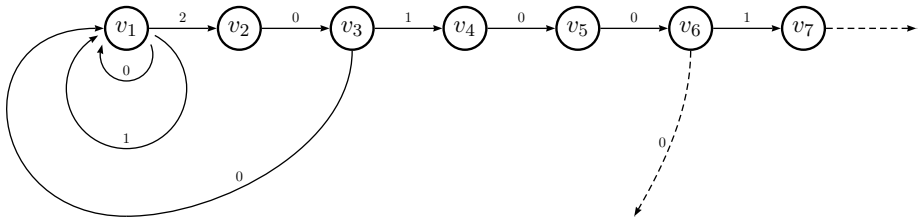

\LargeState\SmallPicture
\VCDraw{%
  \begin{VCPicture}{(-2,-7)(20,2)}
    \State[v_1]{(0,0)}{v1}
    \State[v_2]{(3,0)}{v2}
    \State[v_3]{(6,0)}{v3}
    \State[v_4]{(9,0)}{v4}
    \State[v_5]{(12,0)}{v5}
    \State[v_6]{(15,0)}{v6}
   \State[v_7]{(18,0)}{v7}
    \EdgeL{v1}{v2}{}\taput{2}
    \EdgeL{v2}{v3}{}\taput{0}
    \EdgeL{v3}{v4}{}\taput{1}
    \EdgeL{v4}{v5}{}\taput{0}
   \EdgeL{v5}{v6}{}\taput{0}
    \EdgeL{v6}{v7}{}\taput{1}
    \cnode*(0,-5){0pt}{p0}
    \cnode*(13.5,-5){0}{p1}
    \nccurve[angleA=-90,angleB=60,linestyle=dashed]{->}{v6}{p1}\taput[npos=.7]{0}
    \nccircle[angleA=180,nodesep=5pt]{<-}{v1}{.7cm}\taput[npos=.7]{0}
    \nccircle[angleA=180,nodesep=5pt]{<-}{v1}{1.6cm}\taput[npos=.7]{1}
    \nccurve[angleA=-90,angleB=0]{-}{v3}{p0}\taput[npos=.7]{0}
    \nccurve[angleA=180,angleB=180,ncurv=1.5]{->}{p0}{v1}
    \cnode*(21,0){0pt}{p3}
    \ncline[linestyle=dashed]{v7}{p3}
  \end{VCPicture}}
\label{fig:beta}
\caption{A graph presentation of a $\beta$-shift.}
\end{figure}

Every $\beta$-shift can be presented by a countable state directed labeled graph $\Gamma_\beta$, which is shown in Figure 1 (see also \cite{BH, PfS, Tho7}). We describe the construction of this graph, assuming that $w(\beta)$ is not eventually periodic.  Let $v_1, v_2, \ldots$ be a countable set of vertices.  We use the following two rules to add edges to this graph.  Firstly, for every $i\geq 1$, we draw a directed edge from $v_i$ to $v_{i+1}$ and label it with the value $w_i(\beta)$.  Secondly, if $w_i(\beta) \geq 1$, then for every $j\in \{0, 1, \dots, w_i(\beta) - 1\}$, we draw a directed edge from $v_i$ to $v_1$ labeled with the value $j$.

Note that if $w_i(\beta)=0$, then the only edge which starts at $v_i$ is the edge from $v_i$ to $v_{i+1}$ labeled by $0$, and if $w_i(\beta)>0$ then there is at least one edge from $v_i$ to $v_1$.  The graph $\Gamma_\beta$ characterises $\Sigma_\beta$ as follows:  a sequence $x$ belongs to $\Sigma_\beta$ if and only if $x$ labels an infinite path of directed edges in $\Gamma_\beta$ that starts at the vertex $v_1$.

It follows that words in the language $\LLL$ correspond to finite paths in the graph $\Gamma_\beta$ starting at $v_1$.  Let $\GGG$ be the collection of words for which the corresponding path also ends at $v_1$, and let $\CCC^s$ be the collection of words for which the corresponding path never returns to $v_1$.  Let $\CCC^p = \emptyset$. We can see from the graph that
\[
\CCC^s = \{ w_1(\beta)\cdots w_n(\beta) \mid n \geq 1\}.
\]
This is because the only finite paths that never return to $v_1$ are those which visit the vertices $v_1, v_2, v_3, \dots$ in that order.

Every path can be decomposed into a part that ends at $v_1$ followed by a part that does not return, and so this gives us a decomposition $\LLL = \CCC^p \GGG \CCC^s$.

Loops based at $v_1$ can be freely concatenated and each such loop corresponds to a periodic orbit, hence $\GGG$ has \Fp-specification with $t=0$.  Furthermore, we have
\[
\GGG(M) = \{ vw \mid v\in \GGG, w\in \CCC^s, |w|\leq M \},
\]
and by taking $\tau$ to be such that every path of length $\leq M$ that begins at $v_1$ can return to $v_1$ within $\tau$ steps, we see that Condition~(\ref{extend}) is satisfied.

Furthermore, $\#\CCC^s_n = 1$ for all $n\geq 1$, and so Condition (\ref{zero}$'$) holds.  Thus every $\beta$-shift satisfies the conditions of Theorem~\ref{thm:CPE}, and this proves the part of Theorem~\ref{thm:beta-factor} concerning $\beta$-shifts.

\subsection{$S$-gap shifts}
An $S$-gap shift $\Sigma_S$ is a subshift of $\{0,1\}^\ZZ$ defined by the rule that for a fixed $S \subset \{0,1,2, \ldots\}$, the number of $0$'s between consecutive $1$'s is an integer in $S$. More precisely, the language of $\Sigma_S$ is
\[
\{ 0^n 1 0^{n_1} 1 0^{n_2} 1 \cdots 1 0^{n_k} 1 0^m \mid n_i\in S \text{ for all } 1\leq i\leq k \text{ and } n, m \in \NN \},
\]
together with $\{0^n \mid n \in \NN \}$, where we assume that $S$ is infinite.  (If $S$ is finite, then $\Sigma_S$ is sofic.)  The entropy of the $S$-gap shift is $\log \lambda$, where $\lambda$ is the unique solution to $1 = \sum_{n \in S} x^{-n-1}$.  (See \cite{We} or \cite[Exercise 4.3.7]{LM}.)

The language for $\Sigma_S$ admits the following decomposition:
\begin{align*}
\GGG &= \{ 0^{n_1} 1 0^{n_2} 1 \cdots 1 0^{n_k} 1 \mid n_i\in S \text{ for all } 1\leq i\leq k \}, \\
\CCC^p &= \{ 0^n 1 \mid n\notin S \}, \\
\CCC^s &= \{ 0^n \mid n\in \NN \}.
\end{align*}
It follows immediately that $\GGG$ has \Fp-specification with $t=0$.  Condition (\ref{slow}$'$) follows from the observation that $\#\CCC_n^p \leq 1$ and $\#\CCC_n^s = 1$ for every $n$.

It is not hard to see that Condition \eqref{extend} holds.  (See also the discussion in Section~\ref{examples}.)  Applying Theorem~\ref{thm:CPE} proves the remainder of Theorem~\ref{thm:beta-factor}.

\begin{remark}
If every element of $S$ is odd, then every periodic orbit in the corresponding $S$-gap shift has even period. This demonstrates that in the definition of the periodic orbit measures \eqref{eqn:permeas}, it is crucial that $\Per(n)$ denotes the periodic orbit measures of period at most $n$, rather than those of period exactly $n$.
\end{remark}

\subsection{Specification properties for $S$-gap shifts}
Specification and almost specification properties can be seen to fail for a generic $S$-gap shift. For example, let $S= \{2^n \mid n \in \NN \}$. Consider the words $u = 10^i$ and $w=0^j 1$. Suppose $i+j = 2^n -k > 2^{n-1}$. The shortest word $v$ such that $uvw$ is admissible is $v= 0^k$. Choosing $n$ and $k$ large enough clearly shows that specification fails. Similarly, specification fails for any $S$-gap shift for which $S$ has unbounded gaps (i.e. $\{n_{i+1} - n_i \mid n_i \in S\}$ is unbounded).

We give the definition of almost specification and give a similar example which shows that almost specification fails.

We say a non-decreasing function $g\colon \NN \mapsto \NN$ is a mistake function if $g(n) \leq n$ for all $n$ and $g(n)/n \rightarrow 0$. We say a symbolic space has almost specification if there exists a mistake function $g$ such that for every $w^1, \ldots, w^n \in \LLL$, there exist words $v^1, \ldots, v^n \in \LLL$ with $|v^i| = |w^i|$ such that $v^1 v^2 \cdots v^n \in \LLL$ and each $v^i$ differs from $w^i$ in at most $g(|v^i|)$ places. This is a special case, adapted to symbolic dynamics, of the definition that appears in \cite{PfS, Tho7}.

The space of shifts with almost specification is closed under factors, and every $\beta$-shift has almost specification with the mistake function $g(n)=1$ (see \cite{Tho7} for details). It is an open question whether every shift with almost specification is intrinsically ergodic. We demonstrate that Theorem \ref{thm:main} applies to examples without the almost specification property by giving an example of an $S$-gap shift where almost specification fails.

Consider the $S$-gap shift where $S=\{(2n)! \mid n\in \NN \}$. Let $g(n)$ be any mistake function. Let $n$ be sufficiently large so that $2(g(2n)+1)) \leq 2n$. Let $k= (2n)!$ and let $u = 1^{g(k)+1} 0^{k-g(k)-1}$, and $v= 0^{k-g(k)-1}1^{g(k)+1} $.  Observe that both $u$ and $v$ are in $\LLL$.  Consider the word $uv$. There is no way to change $uv$ into the constant sequence $0^n$. The number of $0$'s we can change to a $1$ is bounded above by $2 g(k)$. Each string of $0$'s we obtain must have a length belonging to $S$. Suppose our modified word has $\ell$ strings of zeroes of length $n_1, \ldots, n_\ell \in S$. Each of these lengths is bounded above by $(2(n-1))!$, and $\ell \leq 2 g(k)+1$. Thus, $\sum_{i=1}^\ell n_i \leq 2(g(k)+1) (2(n-1))! < 2n  (2(n-1))! $. 
There are at most $4g(k)+2 < 8n+2$ entries of $1$. Thus, if we were able to modify $uv$ into an admissible word using $2g(k)$ mistakes, then its length would be strictly less than $2 (2n)!$ which is a contradiction.

\section{Application to coded systems} \label{examples}

A shift space is \emph{coded} if its language $\LLL$ is freely generated by a countable set of words---that is, if there exist words $\{w^n\}_{n\in \NN} \subset \LLL$ such that
\[
\LLL = \overline{\{ w^{n_1} w^{n_2} \cdots w^{n_k} \mid n_i \in \NN \}},
\]
where $\overline{\{\cdot\}}$ denotes closure under the operation of passing to subwords.  We refer to the words $w^n$ as the \emph{generators} of the coded system $(X,\sigma)$. For more information, see \cite{BH, FiFi}.

The language of a coded system has a natural decomposition $\LLL = \CCC^p \GGG \CCC^s$ for which $\GGG$ has \Fp-specification.  Namely, we may consider
\[
\GGG = \{ w^{n_1} \cdots w^{n_k} \mid n_i \in \NN \},
\]
the collection of all concatenations of generators (note that here we do \emph{not} allow passing to subwords), together with
\begin{align*}
\CCC^p &= \{ i_k^s(w^n) \mid n\in \NN, 1\leq k\leq |w^n| \}, \\
\CCC^s &= \{ i_k^p(w^n) \mid n\in \NN, 1\leq k\leq |w^n| \}.
\end{align*}
That is, $\CCC^p$ is the collection of all suffixes of generators, and $\CCC^s$ is the collection of all prefixes of generators.  Conditions (\ref{glue}) and (\ref{extend}) hold for much the same reasons as they did for the $S$-gap shifts.  Indeed, it is immediate that $\GGG$ has \Fp-specification with $t=0$, so Condition \eqref{glue} holds.

For Condition \eqref{extend}, we observe that given $u\in \CCC^p$, there exists a generator $w$ such that $u=i_{|u|}^s(w)$.  Let $\tau^p(u)$ be the minimum value of $|w|$ over all such generators.  Given $u\in \CCC^s$, define $\tau^s(u)$ similarly, as the minimum value of $|w|$ over all generators $w$ such that $u$ is a prefix of $w$.

Given $M\in \NN$, define $\tau^s(M)$ and $\tau^p(M)$ by
\begin{align*}
\tau^s(M) &= \max \{ \tau^s(u) \mid u \in \CCC^p, |u| \leq M \}, \\
\tau^p(M) &= \max \{ \tau^p(u) \mid u \in \CCC^s, |u| \leq M \}.
\end{align*}
Let $\tau(M) = \tau^p(M) + \tau^s(M)$; then given any word $v\in \GGG(M)$, there exist words $x^1, x^2$ with $|x^1| \leq \tau^s(M)$ and $|x^2| \leq \tau^p(M)$ such that $x^1 v x^2 \in \GGG$.  This proves Condition \eqref{extend}.

The number $c_n$ defined before Theorem~\ref{thm:coded} is nothing but $\#(\CCC^p \cup \CCC^s)_n$ and so $\ulim_{n\to\infty} \frac 1n \log c_n = h(\CCC^p \cup \CCC^s)$.  Thus, Theorem~\ref{thm:coded} follows from Theorem~\ref{thm:main} and Theorem~\ref{thm:CPE}. We remark that our results on $\beta$-shifts, $S$-gap shifts, and their factors (Theorem \ref{thm:beta-factor}) are a special case of Theorem~\ref{thm:coded}.

\section{Proof of Theorem~\ref{thm:main}} \label{proof}
\subsection{Uniform estimates on numbers of words}

We obtain estimates on the growth rates of  $\# \LLL_n$ and $\# \GGG_n$.  The estimates in this section require only conditions \eqref{glue} and \eqref{slow}:  condition \eqref{extend} will not be used until we prove the Gibbs property in Section~\ref{sec:gibbs}. The following lemma is a special case of  \cite[Lemma 18.5.3]{KH} or \cite[Lemma 2-3]{Bo7}.

\begin{lemma} \label{lem:Lngeq}
For every $n$,
\begin{equation}\label{eqn:Lngeq}
\#\LLL_n \geq e^{nh(\LLL)}.
\end{equation}
\end{lemma}
\begin{proof}
It is straightforward to obtain $\# \LLL_{m+n} \leq (\#\LLL_m)(\#\LLL_n)$, which yields $\#\LLL_{kn} \leq (\#\LLL_n)^k$, and upon taking logarithms, $\frac 1{kn} \log \#\LLL_{kn} \leq \frac 1n \log \#\LLL_n$.  Passing to the limit as $k\to\infty$ gives the result.
\end{proof}

Using condition (\ref{glue}), we can obtain an upper bound on $\#\GGG_n$.

\begin{lemma} \label{lem:Gnleq}
There exists $C_1 >0$ such that for all $n$,
\begin{equation}\label{eqn:Gnleq}
\#\GGG_n \leq C_1 e^{nh(\LLL)}.
\end{equation}
\end{lemma}
\begin{proof}
Condition (\ref{glue}) immediately implies that $\#\LLL_{k(n+t)} \geq (\#\GGG_n)^k$, which gives
\[
\frac 1{k(n+t)} \log \# \LLL_{k(n+t)} \geq \frac 1{n+t} \#\GGG_n.
\]
Sending $k$ to infinity, we obtain $\#\GGG_n \leq e^{(n+t)h(\LLL)}$.
\end{proof}
This leads to an upper bound on $\#\LLL_n$ by using the decomposition $\LLL = \CCC^p \GGG \CCC^s$ together with the bound on $\CCC^p$ and $\CCC^s$ given in condition (\ref{slow}).

\begin{lemma}\label{lem:Lnleq}
There exists $C_2>0$ such that for every $n$,
\begin{equation}\label{eqn:Lnleq}
\#\LLL_n \leq C_2 e^{nh(\LLL)}.
\end{equation}
\end{lemma}
\begin{proof}
Fix $\eps>0$ such that $h(\CCC^p \cup \CCC^s) < h(\LLL) - \eps$.  Then there exists a constant $C_3$ such that \begin{equation}\label{eqn:CpCs}
\#(\CCC^p_n \cup \CCC^s_n) \leq C_3 e^{n (h(\LLL) - \eps)}
\end{equation}
for every $n$.  For every word $x\in\LLL_n$ there are non-negative integers $i,j,k$ that sum to $n$ and for which $x$ can be decomposed as $uvw$, where $u\in \CCC^p_i$, $v\in \GGG_j$, and $w\in \CCC^s_k$.  Thus
\begin{align*}
\#\LLL_n &\leq \sum_{i+j+k=n} (\#\CCC^p_i)(\#\GGG_j)(\#\CCC^s_k) \\
&\leq C_1 C_3^2 \sum_{i+j+k=n} e^{i(h(\LLL) - \eps)} e^{jh(\LLL)} e^{k(h(\LLL) - \eps)} \\
&= C_1 C_3^2 e^{nh(\LLL)} \sum_{i+j+k=n} e^{-(i+k)\eps} \\
&= C_1 C_3^2 e^{nh(\LLL)} \sum_{m=0}^n \sum_{i=0}^m e^{-m\eps} \\
&\leq C_1 C_3^2 e^{nh(\LLL)} \sum_{m\geq 0} (m+1)e^{-m\eps}.
\end{align*}
The sum converges and is independent of $n$, which completes the proof.
\end{proof}

In place of a lower bound for every $\#\GGG_n$, which is not possible, we obtain the following estimate.

\begin{lemma}\label{lem:good-growth0}
There exist constants $C_4>0$ and $N\in \NN$ such that for every $n\in \NN$, there exists $\ell$ with $n-N \leq \ell \leq n$ such that
\begin{equation}\label{eqn:llimpos0}
\# \GGG_\ell \geq C_4 e^{\ell h(\LLL)}.
\end{equation}
\end{lemma}
\begin{proof}
Given $j\in \NN$, write $a_j = \#\GGG_j e^{-jh(\LLL)}$.  Choose $\eps>0$ and $C_3$ so that~\eqref{eqn:CpCs} holds, as in the proof of Lemma~\ref{lem:Lnleq}.  A similar calculation to the one there gives
\begin{align*}
e^{nh(\LLL)} \leq \# \LLL_n
&\leq \sum_{i+j+k=n} (\#\CCC^p_i) (\# \GGG_j) (\#\CCC^s_k) \\
&\leq C_3^2 \sum_{i+j+k=n} e^{(i+k)(h(\LLL) - \eps)} \# \GGG_j.
\end{align*}
This implies that
\begin{align*} \label{eqn:sum10}
 C_3^{-2} &\leq \sum_{i+j+k=n} e^{-(i+k)\eps} \# \GGG_j e^{-jh(\LLL)} \\ &\leq \sum_{m=0}^n (m+1) e^{-m\eps} a_{n-m}.
\end{align*}
Let $N$ be large enough such that $C_5 := C_3^{-2} - \sum_{m\geq N} (m+1) e^{-m\eps}C_1 > 0$, and let $C_6 = \max \{(m+1) e^{-m\eps} \mid m\in \NN \}$.   Then we have 
\[
C_3^{-2} \leq \left(\sum_{m=0}^{N-1} C_6 a_{n-m}\right) + \left( \sum_{m\geq N} (m+1) e^{-m\eps} C_1 \right),
\]
which yields
\[
C_6 \sum_{m=0}^{N-1} a_{n-m} \geq C_5.
\]
It follows that there exists $n-N \leq \ell \leq n$ such that $a_\ell \geq C_5/(N C_6)$.
\end{proof}

We use a similar argument to obtain the following estimate for $\#\GGG(M)_n$.

\begin{lemma} \label{lem:LnL}
For all $\delta >0$, there exists $M = M(\delta) \in \NN$ such that for all $n$,
\begin{equation}\label{eqn:aegood}
\frac{\# \GGG(M)_n}{\# \LLL_n} \geq 1 - \delta.
\end{equation}
\end{lemma}
\begin{proof}
Let $\eps>0$ and $C_3$ be such that~\eqref{eqn:CpCs} holds, and choose $M$ so that $C_1 C_3^2 \sum_{m>M} (m+1)e^{-m\eps} \leq \delta$.

Once again, we use the decomposition property of $\LLL$ to write any word in $\LLL_n$ as $uvw$, where $u\in \CCC^p_i$, $v\in \GGG_j$, and $w\in \CCC^s_k$, with $i+j+k=n$.  If $i+k\leq M$, then $uvw\in \GGG(M)_n$, and it follows that
\[
\#\LLL_n \leq  \# \GGG(M)_n + \sum_{\substack{ i+j+k=n \\ i+k > M}} (\#\CCC^p_i) (\# \GGG_j) (\#\CCC^s_k).
\]
Now using~\eqref{eqn:Gnleq}, \eqref{eqn:CpCs}, and~\eqref{eqn:Lngeq}, we have
\begin{align*}
\#\LLL_n
&\leq \# \GGG(M)_n + C_1 C_3^2 \sum_{\substack{ i+j+k=n \\ i+k > M}} e^{(i+k)(h(\LLL) - \eps)} e^{jh(\LLL)} \\
&\leq \#\GGG(M)_n + C_1 C_3^2 (\#\LLL_n) \sum_{m > M} (m+1)e^{-m\eps} \\
&\leq \#\GGG(M)_n + \delta (\#\LLL_n).\qedhere
\end{align*}
\end{proof}

An immediate consequence of Lemma~\ref{lem:LnL} is that every collection of words that grows quickly enough has arbitrarily large intersection with $\GGG(M)$, in the following sense.

\begin{lemma}\label{lem:good-growth}
Suppose $\DDD\subset \LLL$ and $C_7>0$ are such that
\begin{equation}\label{eqn:Duniform}
\#\DDD_n \geq C_7 e^{nh(\LLL)}
\end{equation}
for every $n$.  Then there exist constants $M\in \NN$ and $C_8>0$ such that for every $n\in \NN$,
\begin{equation}\label{eqn:llimpos}
\# (\DDD_n \cap \GGG(M)_n) \geq C_8 e^{nh(\LLL)}.
\end{equation}
\end{lemma}
\begin{proof}
Let $\delta >0$ be sufficiently small that $C_8 := C_7 - \delta C_2 > 0$, and let $M=M(\delta)$ be given by Lemma~\ref{lem:LnL}.  Then~\eqref{eqn:Lnleq}, \eqref{eqn:aegood}, and \eqref{eqn:Duniform} yield
\begin{align*}
\#(\DDD_n \cap \GGG(M)_n) &\geq \#\DDD_n - \#(\DDD_n \setminus \GGG(M)_n) \\
&\geq C_7 e^{nh(\LLL)} - \#(\LLL_n \setminus \GGG(M)_n) \\
&\geq C_7 e^{nh(\LLL)} - \delta \#\LLL_n \\
&\geq C_7 e^{nh(\LLL)} - \delta C_2 e^{nh(\LLL)} = C_8 e^{nh(\LLL)}.\qedhere
\end{align*}
\end{proof}

\subsection{Collections of cylinders with uniformly positive measure}

For a collection of words $\DDD$ and a measure $\nu$, we abuse notation slightly and write $\nu(\DDD_n) = \nu\left( \bigcup_{w\in \DDD_n} [w] \right)$, where $[w]$ is the central cylinder defined by $w$.  We also write $\nu(w)$ in place of $\nu([w])$ where it will not cause confusion.

\begin{lemma}\label{lem:posmeas}
For all $\gamma\in (0,1)$ there exists $C_9>0$ such that if $\nu$ is a measure of maximal entropy, $n\in \NN$, and $\DDD_n \subset \LLL_n$ has $\nu(\DDD_n) \geq \gamma$, then
\begin{equation}\label{eqn:posmeas}
\#\DDD_n \geq C_9 e^{nh(\LLL)}.
\end{equation}
\end{lemma}
\begin{proof}
Recall from the definition of measure-theoretic entropy that
\[
h_\nu(\sigma) = \lim_{n\to\infty} \frac 1n H_\nu(\AAA^n) = \inf_{n\geq 1} \frac 1n H_\nu(\AAA^n),
\]
where $\AAA^n$ is the partition of $X$ into $n$-cylinders, and
\[
H_\nu(\AAA^n) = \sum_{w\in \LLL_n} -\nu(w) \log \nu(w).
\]
Because $h_\nu(\sigma) = \htop(X_\LLL,\sigma) = h(\LLL)$, this yields the following inequality for every $n$ (we write $\DDD_n^c$ for the complement of $\DDD_n$ in $\LLL_n$):
\begin{align*}
nh(\LLL) &\leq \sum_{w\in \LLL_n} -\nu(w) \log \nu(w) \\
&= \sum_{w\in \DDD_n} -\nu(w) \log \nu(w) + \sum_{w\in \DDD_n^c} -\nu(w) \log \nu(w).
\end{align*}
Normalising each sum yields
\begin{equation}\label{eqn:nhleq}
\begin{aligned}
nh(\LLL) &\leq \nu(\DDD_n) \left(\sum_{w\in \DDD_n} -\frac{\nu(w)}{\nu(\DDD_n)} \log \left(\frac{\nu(w)}{\nu(\DDD_n)}\right) \right) \\
&\qquad +\nu(\DDD_n^c) \left(\sum_{w\in \DDD_n^c} -\frac{\nu(w)}{\nu(\DDD_n^c)} \log \left(\frac{\nu(w)}{\nu(\DDD_n^c)}\right)\right) \\
&\qquad + (-\nu(\DDD_n)\log \nu(\DDD_n) - \nu(\DDD_n^c) \log \nu(\DDD_n^c)).
\end{aligned}
\end{equation}
Recall that for any non-negative numbers $a_1,\dots a_k$ summing to $1$, we have
\[
\sum_{i=1}^k -a_i \log a_i \leq \log k.
\]
We apply this to the first sum in~\eqref{eqn:nhleq} with the quantities $a_i$ replaced by $\frac{\nu(w)}{\nu(\DDD_n)}$, to the second sum with $a_i$ replaced by $\frac{\nu(w)}{\nu(\DDD_n^c)}$, and to the last line with $a_1=\nu(\DDD_n)$ and $a_2=\nu(\DDD_n^c)$.  This yields
\[
nh(\LLL) \leq \nu(\DDD_n) \log \#\DDD_n + \nu(\DDD_n^c) \log \#(\DDD_n^c) + \log 2,
\]
Lemma~\ref{lem:Lnleq} implies that $\#(\DDD_n^c) \leq \#\LLL_n \leq C_2 e^{nh(\LLL)}$, and so we have
\begin{align*}
nh(\LLL) &\leq \nu(\DDD_n) \log \#\DDD_n + (1-\nu(\DDD_n)) (\log C_2 + nh(\LLL)) + \log 2 \\
&= \nu(\DDD_n) \log \#\DDD_n + \log (2C_2) + nh(\LLL) 
 - \nu(\DDD_n) (\log C_2 + nh(\LLL)).
\end{align*}
Rearranging and using the assumption that $\nu(\DDD_n) \geq \gamma$, this gives
\begin{align*}
\nu(\DDD_n) \log \#\DDD_n &\geq \nu(\DDD_n) (\log C_2 + nh(\LLL)) - \log (2C_2), \\
\log \#\DDD_n &\geq \log C_2 + nh(\LLL) - \frac {\log (2C_2)}{\nu(\DDD_n)} \\
&\geq \log C_2 + nh(\LLL) - \gamma^{-1} \log(2C_2),
\end{align*}
and exponentiating both sides yields~\eqref{eqn:posmeas}.
\end{proof}

\begin{lemma}\label{lem:posmeasgood}
For all $\gamma\in (0,1)$ there exists $C_{10}>0$ and $M\in \NN$ such that if $\nu$ is a measure of maximal entropy, $n\in \NN$, and $\DDD_n \subset \LLL_n$ has $\nu(\DDD_n) \geq \gamma$, then
\begin{equation}\label{eqn:posmeas2}
\#(\DDD_n \cap \GGG(M)) \geq C_{10} e^{nh(\LLL)}.
\end{equation}
\end{lemma}
\begin{proof}
This follows from Lemma~\ref{lem:posmeas} and Lemma~\ref{lem:good-growth}.
\end{proof}

\subsection{A Gibbs property}\label{sec:gibbs}

We build a measure of maximal entropy $\mu$ as a limit of $\delta$-measures $\mu_n$ evenly distributed across $n$-orbits.  In particular, for every $n$ we choose a finite set of points $E_n$ such that for every word $w \in \LLL_n$, the central cylinder $[w]$ contains exactly one element of $E_n$.  Consider the measures $\nu_n$ defined by
\[
\nu_n := \frac 1{\#E_n} \sum_{x\in E_n} \delta_x.
\]
In order to obtain invariant measures, we consider the measures
\begin{equation}\label{eqn:mun}
\mu_n := \frac{1}{n}\sum_{k=0}^{n-1} (\sigma^*)^k \nu_n
\end{equation}
and let $\mu$ be a weak* limit of the sequence $\{\mu_n\}$.
\begin{lemma}
$\mu$ is a measure of maximal entropy.
\end{lemma}
\begin{proof}
This is proved in the second part of ~\cite[Theorem 8.6]{Wa}.
\end{proof}

We prove a Gibbs property for the measure of cylinders corresponding to words in $\GGG$.
\begin{lemma} \label{GibbsG}
There exists $C_{11}>0$ such that for every $n\in \NN$ and $w\in \GGG_n$, we have
\begin{equation}\label{eqn:gibbsG}
\mu([w]) \geq C_{11} e^{-nh(\LLL)}.
\end{equation}
\end{lemma}
\begin{proof}
By Lemma~\ref{lem:good-growth0}, we have $N\in \NN$, $C_4>0$, and a sequence $n_j\nearrow\infty$ with $n_{j+1}-n_j \leq N$ such that for all $j$,
\begin{equation}\label{eqn:gntau'}
\#\GGG_{n_j} \geq  C_4 e^{n_j h(\LLL)}.
\end{equation}
Consider $w \in \GGG_n$. We estimate $\mu_m([w])$ for large $m$ by estimating $\nu_m(\sigma^{-k}([w]))$ first. Let $t \in \NN$ be provided by condition \eqref{glue}.  Fix $k \leq m$. If $k-t\leq N$, let $\ell_1 =0$. Otherwise, let $\ell_1 \in \{n_j\}$ and satisfy $k-t-N \leq \ell_1 \leq k-t$.  Let $\ell_2 \in \{n_j\}$ satisfy $ m-k-t-n-N \leq \ell_2 \leq m-k-t-n$.  If $m-k-t-n < 0$, let $\ell_2=0$.

First assume $\ell_1, \ell_2 >0$. It follows from condition \eqref{glue} that for every $v^1\in \GGG_{\ell_1}$ and $v^2\in \GGG_{\ell_2}$ there exist words $u^1,u^2 \in \LLL$ with $|u^i| = t$ such that $x := v^1 u^1 w u^2 v^2 \in \LLL$.  Extending $x$ by at most $N$ symbols at each end, we obtain a word $y\in \LLL_m$.  Different choices of $v^1 $ and $v^2$ give different words $y$, which shows that
\[
\nu_m(\sigma^{-k}([w])) \geq \frac{(\#\GGG_{\ell_1})(\#\GGG_{\ell_2})}{\#\LLL_m}.
\]
If $\ell_i=0$, this formula still holds by setting $\#\GGG_{\ell_i} = 1$.

Since $\ell_i \in \{n_j\}$, we may use~\eqref{eqn:gntau'} and Lemma \ref{lem:Lnleq} to obtain
\begin{align*}
\nu_m(\sigma^{-k}([w])) &\geq C_4^2C_2^{-1} e^{(\ell_1+\ell_2)h(\LLL)}e^{-m h(\LLL)} \\
&\geq C_4^2C_2^{-1}e^{-2(N+t)h(\LLL)} e^{-n h(\LLL)}.
\end{align*}
Writing $C_{11} = C_4^2C_2^{-1}e^{-2(N+t)h(\LLL)}$ and applying this to~\eqref{eqn:mun} gives
\[
\mu_m([w]) \geq C_{11} e^{-nh(\LLL)},
\]
and passing to the limit gives~\eqref{eqn:gibbsG}.
\end{proof}

Up to this point, we have not used Condition \eqref{extend} at all.  From now on we will use this condition as well, which will allow us to extend the Gibbs property in Lemma~\ref{GibbsG} to cylinders corresponding to words in $\GGG(M)$, with the caveat that the constant in the Gibbs property decays as $M \rightarrow \infty$.

\begin{lemma} \label{Gibbs}
For every $M\in \NN$, there exists a constant $K_M>0$ such that for every $n\in \NN$ and $w\in \GGG(M)_n$, we have
\begin{equation}\label{eqn:gibbs}
\mu([w]) \geq K_M e^{-nh(\LLL)}.
\end{equation}
\end{lemma}
\begin{proof}
Fix $M\in \NN$ and let $\tau$ be given by condition~\eqref{extend}.  Then given $w\in \GGG(M)_n$, there exist words $u, v$ with $|u| \leq \tau, |v| \leq \tau$, so that $uwv \in \GGG$. Since $[uwv] \subset \sigma^{-|u|}[w]$, Lemma \ref{GibbsG} gives
\[
\mu([w]) \geq \mu([uwv]) \geq C_{11} e^{-|uvw| h(\LLL)} \geq C_{11} e^{-2\tau h(\LLL)} e^{-n h(\LLL)}.
\]
Setting $K_M = C_{11} e^{-2\tau h(\LLL)}$ gives \eqref{eqn:gibbs}.
\end{proof}

Finally, we observe that there is a uniform upper bound for the $\mu$-measure of an $n$-cylinder.

\begin{lemma} \label{Gibbsupper}
There exists a constant $C_{12}>0$ such that for every $n\in \NN$ and $w\in \LLL_n$, we have
\begin{equation}\label{eqn:gibblunch}
\mu([w]) \leq C_{12} e^{-nh(\LLL)}.
\end{equation}
\end{lemma}
\begin{proof}
Fix $m >n$ and $k<m-n$. By Lemmas \ref{lem:Lnleq} and \ref{lem:Lngeq}, we have
\begin{align*}
\nu_m(\sigma^{-k}([w])) &\leq \frac{(\#\LLL_k)(\#\LLL_{m-k-n})}{\#\LLL_m} \\
&\leq C_2^2 e^{-nh(\LLL)}
\end{align*}
It follows that $\mu_m([w]) \leq C_2^2 e^{-nh(\LLL)}$. Passing to the limit as $m\to\infty$, we obtain~\eqref{eqn:gibblunch} with $C_{12} = C_2^2$.
\end{proof}

Recall that given a set of words $\DDD_n\subset \LLL_n$, we write $\mu(\DDD_n) = \mu (\bigcup_{w\in \DDD_n} [w])$.

\begin{lemma} \label{Cylapprox}
Let $\delta_1 >0$. There exists $M$ so that for all $n$, any subset $\DDD_n \subset \LLL_n$ satisfies
\[
\mu(\DDD_n\cap \GGG(M)) \geq \mu(\DDD_n) - \delta_1.
\]
\end{lemma}
\begin{proof}
Let $\delta = (C_2 C_{12})^{-1} \delta_1$. Lemma \ref{lem:LnL} provides $M \in \NN$ such that for all $n$, 
\[
\#\GGG(M)_n^c = \#(\LLL_n\setminus \GGG(M)_n) \leq \delta \# \LLL_n.
\]
Combining this with Lemma \ref{Gibbsupper} and Lemma \ref{lem:Lnleq} gives
\[
\mu(\GGG(M)_n^c) \leq \delta  \# \LLL_n C_{12} e^{-nh(\LLL)} \leq \delta C_2 C_{12} = \delta_1,
\]
and we have
\begin{multline*}
\mu(\DDD_n) = \mu(\DDD_n\cap \GGG(M)) + \mu(\DDD_n \cap \GGG(M)^c) \\
\leq \mu(\DDD_n\cap\GGG(M)) + \mu(\GGG(M)_n^c) 
\leq \mu(\DDD_n\cap \GGG(M)) + \delta_1,
\end{multline*}
as required.
\end{proof}

\subsection{Proof that $\mu$ is ergodic}

We need to show that the measure $\mu$ is ergodic.  This is a direct consequence of the following result.

\begin{proposition}\label{prop:reallyweakmix}
If two measurable sets $P, Q \subset X$ both have positive $\mu$-measure, then $\ulim_{n\to\infty} \mu(P \cap \sigma^{-n}(Q)) > 0$.
\end{proposition}
\begin{proof}
We begin by considering the case where $P$ and $Q$ are cylinders corresponding to words in $\GGG$.

\begin{lemma}\label{lem:partialmix}
There exists $C_{13} >0$ and $m_j \rightarrow \infty$ such that if $u, v \in \GGG$, then for all sufficiently large $j$,
\begin{equation}\label{eqn:mixuv}
\mu([u]\cap \sigma^{-m_j}[v]) \geq C_{13} \mu([u]) \mu([v]).
\end{equation}
\end{lemma}
\begin{proof}
As in the proof of Lemma~\ref{GibbsG}, we use Lemma~\ref{lem:good-growth0} to obtain $N\in \NN$, $C_4>0$, and a sequence $n_j\nearrow\infty$ with $n_{j+1}-n_j \leq N$ such that~\eqref{eqn:gntau'} holds for all $j$.  Let $m_j = n_j +2t$.

Consider $u, v \in \GGG$. Let $m \in \NN$ be large and fix $k \leq m$. We estimate
\[
\nu_m(\sigma^{-k}[u] \cap \sigma^{-(k+m_j)} [v])).
\]
By a similar argument to Lemma \ref{GibbsG}, we obtain
\[
\nu_m(\sigma^{-k}[u] \cap \sigma^{-(k+m_j)} [v])) \geq \frac{(\#\GGG_{\ell_1})(\# \GGG_{n_j})(\#\GGG_{\ell_2})}{\#\LLL_m},
\]
where $\ell_1 = n_{i_1}$ for some $i_1$ and satisfies $k-t-N \leq \ell_1 \leq k-t$ (or is $0$ if no such number exists), and similarly, $\ell_2 = n_{i_2}$ and satisfies
\[  m- 3t-m_j- |u|-|v|- N \leq \ell_2 \leq m- 3t-m_j- |u|-|v|,\]
(or $0$ if no such $\ell_2$ exists). Using Lemmas \ref {lem:good-growth0}, \ref{lem:Lnleq}, and \ref{Gibbsupper} we obtain
\begin{align*}
\nu_m(\sigma^{-k}[u] \cap \sigma^{-(k+m_j)} [v]))
&\geq C_4^3 C_2 e^{(\ell_1+\ell_2+n_j)h(\LLL)}e^{-m h(\LLL)} \\
&\geq C_4^3 C_2 e^{-2Nh(\LLL)} e^{-(|u| +|v|) h(\LLL)}\\
&\geq C_4^3 C_2 e^{-2Nh(\LLL)} C_{12}^{-2} \mu([u]) \mu([v]).
\end{align*}
Writing $C_{13} := C_4^3 C_2 e^{-2Nh(\LLL)} C_{12}^{-2}$, this yields
\[
\mu_m([u] \cap \sigma^{-m_j} [v])) \geq C_{13} \mu([u]) \mu([v]),
\]
and passing to the limit as $m\to\infty$ 
gives the required estimate.
\end{proof}

This result immediately extends to unions of cylinders from $\GGG_n$.

\begin{lemma}\label{lem:partialmix1.5}
Let $C_{13}>0$ and $m_j\to\infty$ be as in Lemma~\ref{lem:partialmix}, and consider $P\subset \GGG_n$, $Q\subset \GGG_{n'}$.  Let $[P] = \bigcup_{w\in P} [w]$, and similarly for $Q$.  Then for all sufficiently large $j$,
\begin{equation}\label{eqn:mixPQ}
\mu([P]\cap \sigma^{-m_j}[Q]) \geq C_{13} \mu([P]) \mu([Q]).
\end{equation}
\end{lemma}
\begin{proof}
This is a straightforward computation.
\begin{align*}
\mu([P]\cap \sigma^{-m_j}[Q]) &= \sum_{\substack{w\in P \\ w'\in Q}} \mu([w] \cap \sigma^{-m_j} [w']) \\
&\geq \sum_{\substack{w\in P \\ w'\in Q}} C_{13} \mu([w]) \mu([w'])
= C_{13} \mu([P]) \mu([Q]).\qedhere
\end{align*}
\end{proof}

Using Condition \eqref{extend}, this result generalises to unions of cylinders from $\GGG(M)_n$.  As before, given $P\subset \LLL_n$, we write $\mu(P) = \mu([P]) = \mu(\bigcup_{w\in P} [w])$.

\begin{lemma}\label{lem:partialmix2}
Given $M\in \NN$, there exists a constant $K'_M$ such that for every $P\subset\GGG(M)_n$ and $Q \subset \GGG(M)_{n'}$, we have
\begin{equation}\label{eqn:muPQ2}
\ulim_{m\to\infty} \mu(P \cap \sigma^{-m}(Q)) \geq K'_M \mu(P) \mu(Q).
\end{equation}
\end{lemma}
\begin{proof}
For each $w\in \GGG(M)$, using Condition \eqref{extend}, we can choose $x(w), y(w) \in \LLL$ such that $|x(w)|\leq \tau$, $|y(w)|\leq \tau$, and $x(w)wy(w)\in \GGG$. Given $0\leq i,j,i',j'\leq \tau$, let
\begin{align*}
P(i,j) &= \{w\in P \mid i=|x(w)|, j=|y(w)|\}, \\
Q(i',j') &= \{w'\in Q \mid i'=|x(w')|, j'=|y(w')|\}.
\end{align*}
There exist $i,j,i',j'$ such that
\begin{align*}
\#P(i,j) &\geq (\tau+1)^{-2} \#P, \\
\#Q(i',j') &\geq (\tau+1)^{-2} \#Q.
\end{align*}
Now let
\begin{align*}
\hat P &= \{x(w) w y(w) \mid w\in P(i,j) \} \subset \GGG_{n+i+j}, \\
\hat Q &= \{x(w') w' y(w') \mid w'\in Q(i',j') \} \subset \GGG_{n'+i'+j'}.
\end{align*}
We can estimate the left-hand side of~\eqref{eqn:muPQ2} by observing that
\[
\ulim_{m\to\infty} \mu(P \cap \sigma^{-(m+i'+j)}(Q)) \geq
\ulim_{m\to\infty} \mu(\hat P \cap \sigma^{-m} (\hat Q))
\geq C_{13} \mu(\hat P) \mu(\hat Q),
\]
where the second inequality follows from Lemma~\ref{lem:partialmix1.5}.  Furthermore, the Gibbs properties~\eqref{eqn:gibbsG} and~\eqref{eqn:gibblunch} imply that
\begin{align*}
\mu(\hat P) &\geq \#P(i,j) C_{11} e^{-(n+i+j)h(\LLL)} \\
&\geq (\tau+1)^{-2} \#P C_{11} e^{-(i+j)h(\LLL)} e^{-nh(\LLL)} \\
&\geq (\tau+1)^{-2} C_{11} e^{-2\tau h(\LLL)} C_{12}^{-1} \mu(P),
\end{align*}
A similar estimate on $\mu(\hat Q)$ suffices to complete the proof.
\end{proof}

Lemma~\ref{lem:partialmix2} is the key tool in the proof of the following lemma.

\begin{lemma}\label{lem:kindasortamix}
Suppose $\delta_1>0$ and $M$ are such that Lemma~\ref{Cylapprox} holds and let $K'_M$ be as in Lemma~\ref{lem:partialmix2}.  Then for every pair of measurable sets $P,Q\subset X$, we have
\begin{equation}\label{eqn:kindasortamix}
\ulim_{n\to\infty} \mu(P\cap \sigma^{-n}(Q)) \geq K'_M (\mu(P) - \delta_1)(\mu(Q) - \delta_1).
\end{equation}
\end{lemma}
\begin{proof}
Fix $\eps>0$ and choose sets $U,V$ that are unions of cylinders of the same length and for which $\mu(U\symdiff P) < \eps$ and $\mu(V\symdiff Q)<\eps$.  Let $U'\subset U$ be the union of all cylinders in $U$ corresponding to words in $\GGG(M)$, and similarly for $V'\subset V$.  By Lemma~\ref{Cylapprox}, we have $\mu(U') > \mu(U)-\delta_1$ and $\mu(V')>\mu(V)-\delta_1$, and furthermore, by Lemma~\ref{lem:partialmix2},
\begin{equation}\label{eqn:U'V'}
\ulim_{n\to\infty} \mu(U' \cap \sigma^{-n}(V')) \geq K'_M \mu(U') \mu(V').
\end{equation}
We have $U \cap \sigma^{-n}(V) \supset U' \cap \sigma^{-n}(V')$, and so
\begin{equation}\label{eqn:UV}
\ulim_{n\to\infty} \mu(U \cap \sigma^{-n}(V)) \geq K'_M (\mu(U) - \delta_1)(\mu(V) - \delta_1).
\end{equation}
We also observe that
\begin{align*}
|\mu(U\cap \sigma^{-n}(V)) - \mu(P\cap \sigma^{-n}(Q))| &\leq
\mu((U\cap \sigma^{-n}(V)) \symdiff (P\cap \sigma^{-n}(Q))) \\
&\leq \mu((U\symdiff P) \cap \sigma^{-n}(V \symdiff Q)) < \eps
\end{align*}
for every $n$, which together with~\eqref{eqn:UV} implies
\[
\ulim_{n\to\infty} \mu(P \cap \sigma^{-n}(Q)) \geq K'_M (\mu(P) - \delta_1)(\mu(Q) - \delta_1) - \eps.
\]
Since $\eps>0$ was arbitrary,~\eqref{eqn:kindasortamix} follows.
\end{proof}

Now let $P,Q\subset X$ be any measurable sets with positive $\mu$-measure.  For sufficiently small $\delta_1>0$, the right hand side of~\eqref{eqn:kindasortamix} is positive, which completes the proof of Proposition~\ref{prop:reallyweakmix}.
\end{proof}

\subsection{Contradiction if there is another mme}

Let $\mu$ be the ergodic mme constructed in the previous sections, and suppose that some ergodic measure $\nu \perp \mu$ is such that $h_\nu(\sigma) = \htop(X_\LLL,\sigma) =  h(\LLL)$.  Let $\DDD$ be a collection of words such that $\nu(\DDD_n) \to 1$ and $\mu(\DDD_n) \to 0$.  Applying Lemma~\ref{lem:posmeasgood}, we see that there are constants $C_{10}>0$ and $M\in \NN$ such that
\[
\#(\DDD_n \cap \GGG(M)) \geq C_{10} e^{nh(\LLL)}
\]
for every $n$.  Now we use the Gibbs property~\eqref{eqn:gibbs} to observe that
\[
\mu(\DDD_n) \geq \mu(\DDD_n \cap \GGG(M)) \geq K_M e^{-n h(\LLL)} \# (\DDD_n \cap\GGG(M)) \geq K_M C_{10} > 0,
\]
which contradicts the fact that $\mu(\DDD_n) \to 0$.  This contradiction implies that any mme $\nu$ is absolutely continuous with respect to $\mu$, and since $\mu$ is ergodic, this in turn implies that $\nu=\mu$, which completes the proof of the theorem.

\subsection{Characterisation of the unique mme}
We prove the final statement in Theorem \ref{thm:main} under the assumption that $\GGG$ satisfies \Fp-specification. By Lemma~\ref{lem:good-growth0}, there exist constants $C_4>0$ and $N\in \NN$ such that for every $n$, there exists $n-N \leq \ell \leq n$ for which
\[
\#\GGG_\ell \geq C_4 e^{\ell h(\LLL)}.
\]
By \Fp-specification, every word $w\in \GGG_\ell$ determines a periodic orbit of length $\ell+t\leq n+t$, and so we have
\[
\#\Per(n+t) \geq C_4 e^{\ell h(\LLL)} \geq C_4 e^{(n-N)h(\LLL)}.
\]
This in turn yields
\[
\frac 1{n+t} \log \#\Per(n+t) \geq \frac {\log C_4}{n+t} + \frac{n-N}{n+t} h(\LLL),
\]
and so $\lim_{n\to\infty} \frac 1n \log \#\Per(n) = h(\LLL)$.  Standard arguments such as those in the proof of~\cite[Theorem 8.6]{Wa} show that any limit measure $\nu$ of the sequence $\mu_n$ in~\eqref{eqn:permeas} has $h_\nu(\sigma) = h(\LLL)$. Since we showed that $\mu$ is the unique measure of maximal entropy, this shows that the sequence $\mu_n$ converges to $\mu$.

\section{Proofs of other technical results}\label{factors}

\subsection{Proof of Proposition \ref{prop:measures}}\label{sec:measures}

Let $\epsilon$ denote the empty word, so $\epsilon w=w\epsilon=w$ for every $w\in \LLL$, and $[\epsilon] = X_\LLL$.
A measure $\mu$ induces a function $m\colon \LLL\to [0,\infty)$ by $m(w) = \mu([w])$, and this gives a one-to-one correspondence between $\mc M_\sigma(X_\LLL)$ and functions $m\colon \LLL\to [0,1]$ satisfying
\begin{enumerate}
\item $m(\epsilon)=1$;
\item for every $w\in \LLL$ we have $m(w) = \sum_{a=1}^p m(wa) = \sum_{a=1}^p m(aw)$.
\end{enumerate}

Because, by $\sigma$-invariance, the starting point of the cylinder makes no difference to the measure, there is also a one-to-one correspondence between $\mc M_\sigma(\hat X_\LLL)$ and functions $m$ satisfying the conditions above.  This shows that the invariant measures of $X_\LLL$ and $\hat X_\LLL$ can be identified.  Furthermore, the entropy of $\mu$ is determined by $m$, so this identification preserves entropy.

\subsection{Proof of Proposition \ref{prop:factorcat}}

We study the behaviour of the decomposition $\LLL = \CCC^p \GGG \CCC^s$ and its properties under factors, proving Proposition \ref{prop:factorcat}.

In the following, we sometimes write $v\cdot w$ in place of $vw$ to denote concatenation.  Let $\Sigma \subset \Sigma_p$ and $\tilde\Sigma \subset \Sigma_{\tilde p}$ be arbitrary closed two-sided invariant subshifts (the one-sided case is similar), and suppose that $\tilde\Sigma$ is a topological factor of $\Sigma$---that is, there exists a continuous and surjective map $\pi\colon \Sigma \mapsto \tilde\Sigma$ such that $\sigma \circ \pi = \pi \circ \sigma$.  By the fundamental result of Curtis--Lyndon--Hedlund~\cite[Theorem 6.29]{LM}, $\pi$ is a block code: there exist $k\in \NN$ and $\phi\colon \LLL_{2k+1} \to \{1, \ldots, \tilde p \}$ such that
\[
(\pi x)_n = \phi(x_{n-k} x_{n-k+1} \cdots x_{n+k-1} x_{n+k}).
\]
This induces a map $\Phi\colon \LLL_{n+2k} \to \tilde\LLL_n$ by
\[
\Phi(w_1 \cdots w_{n+2k}) = \phi(w_1 \cdots w_{2k+1}) \phi(w_2 \cdots w_{2k+2}) \cdots \phi(w_n \cdots w_{n+2k}).
\]
The map $\Phi\colon \LLL\to \tilde\LLL$ has the following important properties:
\begin{enumerate}
\item $\Phi$ is surjective.
\item For every word $w\in \LLL$ and $x\in {}_n[w]$, we have $\pi(x) \in {}_{n+k}[\Phi(w)]$.
\end{enumerate}
If $k=0$, then $\Phi$ is a homomorphism in the sense that $\Phi(vw)=\Phi(v)\Phi(w)$ for all words $v,w\in \LLL$.  For $k>0$, we need to define maps on $\LLL$ that extract prefixes and suffixes:  in particularly, consider maps $i_k^p, i_k^s\colon \LLL_{\geq k} \to \LLL_k$ given by
\begin{align*}
i_k^p(w) &= w_1 \dots w_k, \\
i_k^s(w) &= w_{|w|-k+1} \dots w_{|w|}.
\end{align*}
Now $\Phi$ has the property that for every $v,w\in \LLL$, we have
\begin{equation}\label{eqn:Phi}
\begin{aligned}
\Phi(vw) &= \Phi(v) \Phi(i_{2k}^s(v) \cdot i_{2k}^p(w)) \Phi(w) \\
&= \Phi(v) \Phi(i_{2k}^s(v) \cdot w) \\
&= \Phi(v \cdot i_{2k}^p(w)) \Phi(w).
\end{aligned}
\end{equation}

Given a decomposition $\LLL = \CCC^p \GGG \CCC^s$, the obvious thing to do is to define a decomposition of $\tilde\LLL$ by applying $\Phi$ to each of $\CCC^p$, $\GGG$, and $\CCC^s$.  Because the homomorphism property of $\Phi$ takes the form~\eqref{eqn:Phi}, we must alter this slightly and define subsets of $\tilde\LLL$ by
\begin{align*}
\tilde\GGG &= \Phi(\GGG), \\
\tilde\CCC^p &= \Phi(\CCC^p \cdot i_{2k}^p(\GGG)), \\
\tilde\CCC^s &= \Phi(i_{2k}^s(\GGG) \cdot \CCC^s).
\end{align*}
Given $u\in \CCC^p$, $v\in \GGG$, and $w\in \CCC^s$, we see from~\eqref{eqn:Phi} that
\[
\Phi(uvw) = \Phi(uv) \Phi(i_{2k}^s(v) \cdot w)
= \Phi(u \cdot i_{2k}^p(v)) \Phi(v) \Phi(i_{2k}^s(v) \cdot w) \in \tilde\CCC^p \tilde\GGG \tilde\CCC^s.
\]
This gives the decomposition of $\tilde \LLL$ claimed in Proposition~\ref{prop:factorcat}.  Furthermore, we observe that for every $M$, we have
\begin{align*}
\Phi(\GGG(M)) &= \{ \Phi(uvw) \mid u\in \CCC^p, v\in \GGG, w\in \CCC^s, |u|\leq M, |w|\leq M \} \\
&= \{ \Phi(u \cdot i_{2k}^p(v)) \Phi(v) \Phi(i_{2k}^s(v) \cdot w) \mid \\
&\qquad\qquad\qquad u\in \CCC^p, v\in \GGG, w\in \CCC^s,
|u|\leq M, |w|\leq M \} \\
&= \tilde\GGG(M).
\end{align*}

Now suppose $\GGG$ has \F-specification.  Then given $w^1,\ldots,w^n\in \tilde\GGG$, we have $w^j = \Phi(v^j)$ for some $v^j \in \GGG$, and by (\ref{glue}) there exist $x^1,\ldots,x^{n-1} \in \LLL_t$ such that
\[
v^1 x^1 v^2 x^2 \cdots x^{n-1} v^n \in \LLL.
\]
Applying $\Phi$ and writing $y^j = \Phi( i_{2k}^s(v^j) \cdot x^j \cdot i_{2k}^p(v^{j+1}))$, we have
\[
w^1 y^1 w^2 y^2 \cdots y^{n-1} w^n \in \LLL.
\]
Thus $\tilde\GGG$ has \F-specification with connecting words of length $t+2k$.  Since $\pi$ takes periodic orbits to periodic orbits, \Fp-specification is also preserved by $\Phi$.

Condition \eqref{extend} is clearly preserved by $\Phi$. If $\GGG(M)$ satisfies Condition \eqref{extend} then for every $v \in \tilde\GGG(M)$, there exists $u, v$ with $|u| \leq t+2k$ and $|w| \leq t+2k$ such that $uvw \in \tilde \GGG$.

Finally, we observe that $\#\tilde\CCC^p_n \leq (\#\LLL_{2k})(\#\CCC^p_n)$, and similarly for $\CCC^s$, which completes the proof of Proposition~\ref{prop:factorcat}.

\subsection{Proof of Proposition~\ref{prop:posent}}

We show that \F-specification implies positive entropy unless the collection $\GGG$ has a very specific structure.

\begin{lemma}\label{lem:vw=wv}
Let $(X,\sigma)$ be a shift space whose language $\LLL$ contains a collection of words $\GGG\subset \LLL$ with \F-specification, and suppose that $\htop(X,\sigma)=0$.  Then given any $v,w\in \GGG$, we have $v \LLL_t w \cap w\LLL_t v \neq \emptyset$; that is, there exist $y,z\in \LLL_t$ such that $vyw=wzv$.
\end{lemma}
\begin{proof}
Fix $v\neq w\in \GGG$, and suppose $v\LLL_t w \cap w\LLL_t v =\emptyset$.  We will conclude that $\htop(X,\sigma)>0$, which suffices to prove the lemma.

Given $N\in \NN$, consider $\xi \in \{1,2\}^N$.  Using \F-specification, there exists $\alpha(\xi) \in \LLL$ having the form
\[
\alpha(\xi) = a^1 b^1 a^2 b^2 \cdots b^{N-1} a^N,
\]
where $b^i\in \LLL_t$, and where
\[
a^i \in \begin{cases}
v \LLL_t w & i=1, \\
w \LLL_t v & i=2.
\end{cases}
\]
Because $v\LLL_t w \cap w\LLL_t v=\emptyset$, the map $\alpha$ is injective.  Furthermore, writing $m=|v|$ and $n=|w|$, each word $\alpha(\xi)$ has length $N(m+n+2t) - t$, and it follows that
\[
\frac 1{N(m+n+2t) - t} \log \#\LLL_{N(m+n+2t) - t} \geq \frac 1{N(m+n+2t) - t} \log 2^N.
\]
Taking a limit as $N\to\infty$ yields $\htop(X,\sigma)\geq \frac 1{m+n+2t} \log 2 > 0$, which contradicts the zero entropy assumption.
\end{proof}

From now on we assume that $\htop(X,\sigma)=0$, and aim to show that $X$ comprises a single periodic orbit.  First we show that every word in $\GGG$ is a prefix of a single infinite sequence.

\begin{lemma}\label{lem:allprefix}
There exists $\hat x\in X$ such that every word $w\in \GGG$ is of the form $w=\hat x_1\cdots \hat x_n$ for some $n$.
\end{lemma}
\begin{proof}
It suffices to show that for any $v,w\in \GGG$ with $m=|v|\leq |w|$, we have $w_1\cdots w_m=v$.  This follows from Lemma~\ref{lem:vw=wv}, since there are words $y,z\in \LLL_t$ such that $vyw=wzv$, and comparing the first $m$ symbols of this common word gives the result.
\end{proof}

\begin{lemma}\label{lem:partperiod}
Given $v,w\in \GGG$, write $m=|v|$, $n = |w|$, and suppose that $n > m$.  Then $w$ is $(m+t)$-periodic; that is, $w_i = w_j$ whenever $j\equiv i \bmod (m+t)$.
\end{lemma}
\begin{proof}
We show that $w$ is $(k+t)$-periodic, where $k+t = \gcd(m+t,n+t)$.  Let $u^i \in \LLL_k$ and $x^i\in \LLL_t$ be such that $\hat x = u^1 x^1 u^2 x^2 \cdots$, where $\hat x$ is the sequence from Lemma~\ref{lem:allprefix}.  Write $a=\frac{m+t}{k+t}$ and $b= \frac{n+t}{k+t}$; then $v=u^1x^1u^2 \cdots x^{a-1}u^a$ and $w=u^1x^1u^2 \cdots x^{b-1} u^b$.

By Lemma~\ref{lem:vw=wv}, there are words $y,z\in \LLL_t$ such that $vyw=wzv$. Therefore, the following two expressions represent the same word:
\begin{align*}
&u^1 x^1 \cdots x^{a-1} u^a \,y\ \ u^1\ \ x^1\ \cdots x^{b-a-1} u^{b-a} x^{b-a} u^{b-a+1} \cdots u^b, \\
&u^1 x^1 \cdots x^{a-1} u^a x^a u^{a+1} x^{a+1} \cdots \ x^{b-1} \ \ u^b \ \ z \ \  \ \ u^1 \ \ \ \cdots \ u^a.
\end{align*}
Comparing the subwords $u^i$ and $x^i$ in these two expressions, we see that $u^i = u^j$ and $x^i = x^j$ whenever $i\equiv j\bmod a$.  This follows from a comparison of the two middle segments, where the words $w$ overlap.

Similarly, a comparison of the final segments, where $v$ appears as a suffix of $w$, shows that $u^{b-i} = u^{a-i}$ for all $0\leq i < a$ and $x^{b-i} = x^{a-i}$ for all $1\leq i < a$.  Since $\gcd(a,b)=1$, this is enough to show that $u^1=u^2=\cdots u^b$ and $x^1 = \cdots = x^{b-1}$.
\end{proof}

We use Lemma~\ref{lem:partperiod} to show that the sequence $\hat x\in X$ constructed in Lemma~\ref{lem:allprefix} is periodic.  Indeed, if we fix $v\in \GGG$ and let $m=|v|$, then $\hat x$ is $(m+t)$-periodic.  To see this, observe that Lemma~\ref{lem:partperiod} establishes $(m+t)$-periodicity for every $w\in \GGG$ with $|w| > m$.  By Condition~\eqref{extend}, $\GGG$ contains arbitrarily long words; consequently, there are arbitrarily large values of $n$ such that $\hat x_1 \cdots \hat x_n$ is $(m+t)$-periodic, and this completes the proof.

\section*{Acknowledgements}
We would like to thank Mike Boyle for suggesting this problem, and for making numerous helpful suggestions to improve an early version of this manuscript. We would also like to thank the referee for a careful reading of our manuscript and some helpful suggestions.

\bibliographystyle{amsalpha}
\bibliography{master}

\end{document}